\documentclass[12pt]{article}
\usepackage[utf8]{inputenc}
\usepackage[T1]{fontenc}
\usepackage{fullpage}
\usepackage{amsmath, amssymb, amsthm}

\usepackage{indentfirst}
\usepackage[numbers]{natbib}

\newtheorem{theorem}{Theorem}
\newtheorem{lemma}{Lemma}[section]

\newtheorem{Claim}{Claim}

\usepackage{hyperref}
\usepackage{hyperref}
\usepackage{xcolor}
\usepackage{tikz}
\hypersetup{
colorlinks=true,
linkcolor=blue,
citecolor=blue,
urlcolor=blue,
breaklinks=true
}

\usepackage{etoolbox}
\newtheorem{Subclaim}{Subclaim}
\newcounter{claimsnapshot}

\AtBeginEnvironment{Claim}{\setcounter{Subclaim}{0}}
\AtEndEnvironment{Claim}{\setcounter{claimsnapshot}{\value{Claim}}}
\AtBeginEnvironment{claim}{\setcounter{Subclaim}{0}}
\AtEndEnvironment{claim}{\setcounter{claimsnapshot}{\value{Claim}}}

\begin{document}

\begin{center}	
{\Large \bf When do local bounds permit optimal edge-colorings?}
\vspace{6mm}
		
Xuli Qi*,\,\, Minli Yu

\vspace{3mm}
\baselineskip=0.25in
{\it Department of Mathematics, Hebei University of Science and Technology,}  \\
{\it Shijiazhuang 050018, P.R. China}\\
\vspace{3mm}
\end{center}	
\footnotetext{E-mail address: xuliqi@hebust.edu.cn}

\begin{abstract}
Let $G$ be a multigraph with maximum degree $\Delta(G)$, maximum  multiplicity $\mu(G)$, and chromatic index $\chi'(G)$. A proper edge-coloring is {\bf locally bounded} if each edge $e$ in $G$ with endvertices $u$ and $v$ is assigned a color at most $\max\{d_G(u)+\mu_G(u),d_G(v)+\mu_G(v)\}$, where $d_G(w)$ and $\mu_G(w)$ respectively denote the degree of $w$ and the maximum multiplicity of an edge incident with $w$. The well-known Vizing's theorem guarantees a proper edge-coloring with at most $\Delta(G)+\mu(G)$ colors.
Recently, the local version of this theorem was established for simple graphs by Christiansen~\cite{Christiansen2023}, resolving a conjecture posed in~\cite{BDLP2024}, and extended to multigraphs by Conley, Greb\'ik, and Pikhurko~\cite{CGP2024}. 
In this paper, we define  the  {\bf local chromatic index} of $G$, denoted by $\chi'_{loc}(G)$, as the least positive integer $k$ for which $G$ admits a locally bounded proper edge-coloring with colors in $\{1,2,...,k\}$.
Thus $\Delta(G)\le\chi'(G)\le\chi'_{loc}(G)\le\Delta(G)+\mu(G)$.
We give equivalent characterizations and sufficient conditions for the equality $\chi'_{loc}(G)=\chi'(G)$, identifying when the local bounds in Vizing's theorem can be satisfied using the minimum possible number of colors. In particular, we prove this equality for every regular multigraph, every multigraph $G$ with $\Delta(G)\le 4$, and every simple graph in which the vertices of degree at least five induce a subgraph of maximum degree at most two. Moreover, we construct a simple graph $G$ with $\chi'(G)=\Delta(G)=5$ and $\chi'_{loc}(G)=6$, showing that these local bounds can force the use of an additional color even in simple graphs. Consequently, five is the smallest possible maximum degree for which satisfying these local bounds can require more than $\chi'(G)$ colors.

\noindent\textbf{Keywords:} proper edge-coloring, locally bounded edge-coloring, local chromatic index, local version of Vizing’s theorem
\end{abstract}

\section{Introduction}	
Throughout this paper, all multigraphs are finite, connected, and loopless, and have at least one edge. Multiple edges are allowed unless a multigraph is explicitly assumed to be simple. The standing assumptions of connectedness and a nonempty edge set apply to the graph under study, not to its auxiliary subgraphs or edge-deleted graphs. Let $G=(V(G),E(G))$ be a multigraph, where $V(G)$ is the vertex set and $E(G)$ is the edge set of $G$.  For any positive integer $k$, let $[k]:=\{1,2,\dots,k\}$, that is, the set of the first $k$ positive integers. A {\bf proper $k$-edge-coloring} of a multigraph $G$ is a map $\varphi:E(G)\to[k]$ that assigns to every edge $e$ of $G$ a color $\varphi(e)\in[k]$ such that no two adjacent edges of $G$ receive the same color.
Denote by $\mathcal{C}^k(G)$ the set of all proper $k$-edge-colorings of $G$. The {\bf chromatic index} $\chi'(G)$ is the least integer $k$ such that $\mathcal{C}^k(G)\ne\emptyset$.

Let $\Delta(G)$ denote the maximum degree of $G$, and we  simply set $\Delta:=\Delta(G)$ in this paper. The core of $G$, denoted by $G_\Delta$, is the subgraph of $G$ induced by the set of vertices with the maximum degree in $G$. A graph $G$ is regular if all vertices of $G$ have the same degree.
For $u,v\in V(G)$, we define $E_G(u,v)$ to be the set of edges between $u$ and $v$. In particular, for a simple graph, we identify an edge with its unordered pair of endvertices and write $e=uv$. Furthermore, we define the multiplicity of $\{u,v\}$ to be $\mu_G(u,v):=|E_G(u,v)|$, the maximum edge multiplicity at $u$ to be $\mu_G(u):=\max_{v\neq u}\mu_G(u,v)$, and the maximum multiplicity  of $G$ to be $\mu(G):=\max_{u\neq v}\mu_G(u,v)$. The degree of a vertex $u$ in $G$ is $d_G(u):=\sum_{v\neq u} \mu_G(u,v)$.
	
In the influential paper by Erdős et al.~\cite{Erdos1980}, some local versions of edge-coloring problems were introduced, in which the set of admissible colors for an edge $e\in E_G(u,v)$ is determined by local information such as the degrees of $u$ and $v$, instead of the global maximum degree. Vizing's theorem~\cite{Gupta1967,Vizing1964} states that every graph $G$ admits a proper edge-coloring with at most $\Delta(G)+\mu(G)$ colors. Recently, a local version of this theorem, due to Christiansen~\cite{Christiansen2023} for simple graphs resolving a conjecture posed in~\cite{BDLP2024} and subsequently extended to multigraphs by Conley, Greb\'ik, and Pikhurko~\cite{CGP2024}, guarantees a proper ($\Delta(G)+\mu(G)$)-edge-coloring $\varphi$ satisfying, for every $e\in E(G)$ with endvertices $u$ and $v$,
\begin{equation}
\varphi(e)\le\max\{d_G(u)+\mu_G(u),d_G(v)+\mu_G(v)\}.
\label{eq:local-bound}
\end{equation}
Let $b_G(e):=\max\{d_G(u)+\mu_G(u),d_G(v)+\mu_G(v)\}$ for $e\in E_G(u,v)$ in $G$. We refer to an edge-coloring that satisfies~\eqref{eq:local-bound} as {\bf locally bounded}, and we define the {\bf local chromatic index} of $G$, denoted by $\chi'_{loc}(G)$, as the least positive integer $k$ for which $G$ admits a locally bounded proper $k$-edge-coloring. The local version of Vizing's theorem ensures that this parameter is well defined. Furthermore, we have
\begin{equation}
\Delta(G)\le\chi'(G)\le\chi'_{loc}(G)\le\Delta(G)+\mu(G).
\label{eq:basic-bounds}
\end{equation}
In particular, if $G$ is a simple graph, then $\Delta(G)\le\chi'(G)\leq\chi'_{loc}(G)\leq\Delta(G)+1$.

Dhawan~\cite{Dhawan2025} established list-coloring extensions based on intersections of incident edge lists, yielding local analogues of Vizing’s, Shannon’s, and König’s theorems. We ask when the prescribed numerical local bounds can be satisfied using the optimal palette $[\chi'(G)]$. In this paper, we investigate conditions under which $\chi'_{loc}(G)=\chi'(G)$. Theorem $1$ gives several equivalent characterizations of this equality, while Theorems $2$ and $3$ provide several sufficient conditions.

\begin{theorem} \label{T1}
Let $G$ be a multigraph with $\chi'(G)=k$. The following three conditions are equivalent:
\begin{enumerate}
\item[\upshape(1)] $\chi'_{loc}(G)=\chi'(G)=k$;
\item[\upshape(2)] There exists a proper $k$-edge-coloring $\varphi$ such that $|\{\varphi(e):e\in E(G[D_t])\}|\le t$ for every $t\in[k]$, where $D_t: = \{u\in V(G):d_G(u)+\mu_G(u)\le t\}$;
\item[\upshape(3)] $E(G)$ can be decomposed into $k$ nonempty matchings, i.e., $E(G)=M_1\dot{\cup}\dots\dot{\cup}M_k$, with $b(M_i):=\min_{e\in M_i}b_G(e)$ for $i\in[k]$. After relabeling the matchings so that the values $b_i = b(M_i)$ satisfy $b_1\le b_2 \le\dots\le b_k$, we have $b_i\ge i$ for all $i\in[k]$.
\end{enumerate}
\end{theorem}
	
Borodin et al.~\cite{BKW97} proved that if $G$ is a bipartite multigraph, then there exists a proper edge-coloring $\varphi$ such that $\varphi(e)\le\max\{d_G(u),d_G(v)\}$ for every edge $e\in E_G(u,v)$. Consequently, it follows that $\Delta(G)\le\chi'(G)\le\chi'_{loc}(G)\le \Delta(G)$, which immediately implies $\chi'_{loc}(G)=\chi'(G)$. We prove the following sufficient conditions.
	
\begin{theorem} \label{T2}
Let $G$ be a multigraph with $\chi'(G)=k$.
	
\begin{enumerate}
\item[\upshape(1)]If $G$ is regular, then $\chi'_{loc}(G)=\chi'(G)$;	
\item[\upshape(2)]If $\Delta(G)\le4$, then $\chi'_{loc}(G)=\chi'(G)$;
\item[\upshape(3)]If $d_G(u)+d_G(v)\le b_G(e)+\mu_G(u,v)$ for every edge $e\in E_G(u,v)$ satisfying $b_G(e)<k$, then $\chi'_{loc}(G)=\chi'(G)$.
			
\end{enumerate}
\end{theorem}
	
For simple graphs $G$, Conley, Greb\'ik, and Pikhurko~\cite{CGP2024} showed that if the core $G_\Delta$ has no cycles, then $\chi'_{loc}(G)=\Delta(G)$.
Recall that a simple graph $G$ is of {\bf Class~1 }if $\chi'(G)=\Delta(G)$, and of {\bf Class~2 }if $\chi'(G)=\Delta(G)+1$ \cite{FW1977}. For a Class~2 graph, \eqref{eq:basic-bounds} immediately gives $\chi'_{loc}(G)=\chi'(G)$. Thus, the question is whether a Class~1 graph admits a locally bounded proper edge-coloring with only $\Delta(G)$ colors.
We address this question under a structural restriction on the subgraph induced by vertices of degree at least five. For simple graphs $G$ with maximum degree $\Delta$, let $V_{\ge t}(G):=\{u\in V(G):d_G(u)\ge t\}$, where $t$ is a positive integer.
The following theorem establishes the equality $\chi'_{loc}(G)=\chi'(G)$ for simple graphs whenever every vertex in
$V_{\ge5}(G)$ has at most two neighbors in $V_{\ge5}(G)$.
	
\begin{theorem}  \label{T3}
Let $G$ be a simple graph. If $ \Delta\bigl(G[V_{\ge5}(G)]\bigr)\le2 $,
then $\chi'_{loc}(G)=\chi'(G)$.
\end{theorem}
Equivalently, under the hypothesis of Theorem~3, there exists a proper edge-coloring $\varphi:E(G)\rightarrow[\chi'(G)]$ such that $\varphi(uv)\le\max\{d_G(u),d_G(v)\}+1 $ for every $uv\in E(G)$.
However, the equality $\chi'_{loc}(G)=\chi'(G)$ need not hold for all simple graphs. We construct a counterexample $G^*$ to show that $6=\chi'_{loc}(G^*)>\chi'(G^*)=\Delta(G^*)=5$.

\begin{figure}[htbp]
\begin{center}
\centering
\scalebox{0.55}{\includegraphics{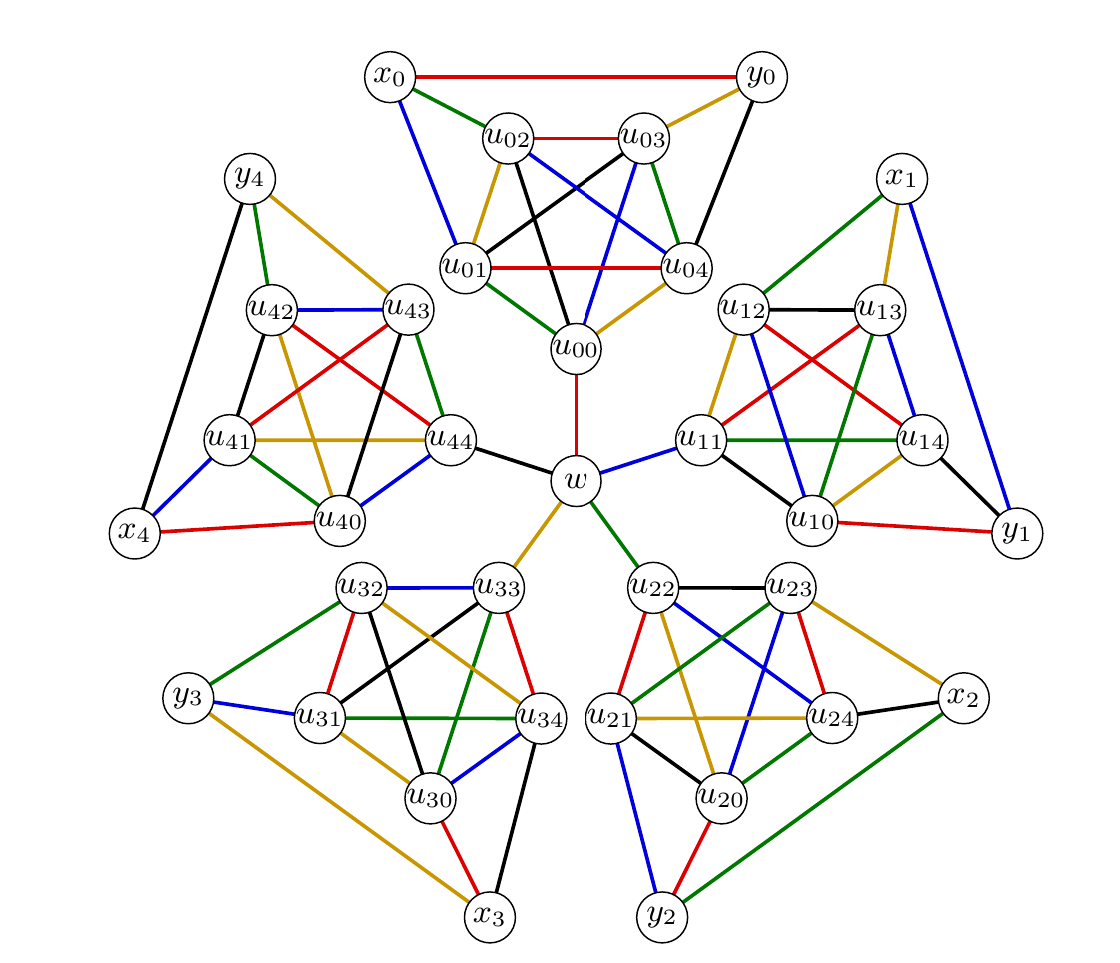}}
\caption{The graph $G^*$ in Example $1$.}
\label{fig:graph}
\end{center}
\end{figure}

\noindent \textbf{Example 1.} Consider the simple graph $G^*$ depicted in Figure 1. Let $V(G^*)=\{w\}\cup\{x_i,y_i,u_{ij}\colon i,j\in\{0,1,2,3,4\}\}$. For each $i\in\{0,1,2,3,4\}$, any two distinct vertices in $\{u_{i0},\dots,u_{i4}\}$ are joined by an edge. In addition, include the edges $wu_{ii}$, $x_iu_{i(i+1)}$, $x_iu_{i(i+2)}$, $y_iu_{i(i+3)}$, $y_iu_{i(i+4)}$, and $x_iy_i$, where all subscripts are computed modulo 5. We claim that $6=\chi'_{loc}(G^*)>\chi'(G^*)=\Delta(G^*)=5$.

It is easily checked that $\chi'(G^*)=\Delta(G^*)= 5$ (see Figure 1). Let $\varphi$ be an arbitrary proper 5-edge-coloring of $G^*$. For each $i\in\{0,1,2,3,4\}$, let $g_{ii}=wu_{ii}$, $g_{i(i+1)}=x_iu_{i(i+1)}$, $g_{i(i+2)}=x_iu_{i(i+2)}$, $g_{i(i+3)}=y_iu_{i(i+3)}$, $g_{i(i+4)}=y_iu_{i(i+4)}$, and $f_i=x_iy_i$. We claim that $g_{ii},\dots,g_{i(i+4)}$ receive pairwise distinct colors and $\varphi(f_i)=\varphi(g_{ii})$ for each $i\in\{0,1,2,3,4\}$.
Fix $i\in\{0,1,2,3,4\}$. For each color $\alpha\in[5]$, let $n_\alpha =|\{j:\varphi(g_{ij})=\alpha\}|$ and $m_\alpha =|\{u_{is}u_{it}:\varphi(u_{is}u_{it})=\alpha\text{ with }0\le s<t\le 4\}|$. For each vertex in $\{u_{i0},\ldots,u_{i4}\}$, its degree is $5$ and five incident edges use exactly all five colors in $[5]$. Each vertex of $\{u_{i0},\ldots,u_{i4}\}$ is incident with an edge colored with $\alpha$. Each edge $u_{is}u_{it}$ colored with $\alpha$ covers two vertices, and each edge $g_{ij}$ colored with $\alpha$ covers one vertex in $\{u_{i0},\ldots,u_{i4}\}$. Hence $2m_\alpha+n_\alpha=5$. Thus $n_\alpha$ is a positive odd integer, and $n_\alpha \ge 1$. Since the five edges $g_{ij}$ satisfy $\sum_{\alpha=1}^5 n_\alpha = 5$, it follows that $n_\alpha=1$ for all $\alpha\in[5]$. Therefore, we obtain that $g_{ii},\dots,g_{i(i+4)}$ receive pairwise distinct colors from $[5]$. Note that the edges adjacent to $f_i$ are exactly $g_{i(i+1)},g_{i(i+2)},g_{i(i+3)}$, and $g_{i(i+4)}$. Thus we obtain $\varphi(f_i)=\varphi(g_{ii})$.

Since $d_{G^*}(w)=5$, $g_{00},\dots,g_{44}$ receive distinct colors from $[5]$. Thus there must exist an edge $f_i$ colored with color 5, where $i\in\{0,1,2,3,4\}$. Note that $d_{G^*}(x_i)=d_{G^*}(y_i)= 3$. We obtain $\varphi(f_i)=5>\max\{d_{G^*}(x_i),d_{G^*}(y_i)\}+1=4$. Therefore, $G^*$ admits no locally bounded proper 5-edge-coloring, and so $\chi'_{loc}(G^*) \ge 6$. Since $G^* $ is a simple graph and $\Delta(G^*)=5$, we have $\chi'_{loc}(G^*)=6>5=\chi'(G^*)$ as desired.

\section{Definitions and preliminary results}
	
A multigraph $G$ is {\bf edge-chromatic critical} if $\chi'(H)<\chi'(G)$ for every proper subgraph $H$ of $G$. Furthermore, we call $G$ a {\bf $\Delta$-critical} graph if $G$ is edge-chromatic critical and $\chi'(G)=\Delta(G)+1$, where $\Delta:=\Delta(G)$. For $e \in E(G)$, $G-e$ represents the graph obtained from $G$ by deleting the edge $e$. An edge $e$ is a {\bf critical edge} of $G$ if $\chi'(G-e)<\chi'(G)$. Clearly, if $G$ has no isolated vertex and every edge is critical, then $G$ is edge-chromatic critical. For a proper $k$-edge-coloring $\varphi$, an edge $e\in E_G(u,v)$ is defined as a {\bf non-local $i$-edge} if $\varphi(e)=i>\max\{d_G(u)+\mu_G(u), d_G(v)+\mu_G(v)\}$, where $i\in[k]$. For every $i\in[k]$, let $B^\varphi_i:=\{e:\varphi(e)=i>b_G(e)\}$ denote the set of all non-local $i$-edges. In particular, if $G$ is a simple graph, then $B^\varphi_i:=\{uv:\varphi(uv)=i,\max\{d_G(u),d_G(v)\}\leq i-2\}$. An {\bf $i$-vertex} in $G$ is a vertex with degree $i$. For $u \in V(G)$, an {\bf $i$-neighbor} of $u$ is a neighbor of $u$ that is an $i$-vertex in $G$. Let $N_G(u)$ denote the set of neighbors of $u$ in $G$. For $u,v\in V(G)$, the distance between $u$ and $v$ is defined as the number of edges on a shortest path between $u$ and $v$ in $G$. If a multigraph $G$ satisfies $\left|E(G)\right|>\Delta(G)\lfloor \frac{|V(G)|}{2}\rfloor$, then $\chi'(G)\ge\Delta(G)+1$. Such graphs are called {\bf overfull}.
	
Let $G$ be a multigraph with a proper $k$-edge-coloring $\varphi\in \mathcal{C}^k(G)$ for some integer $k\ge\Delta(G)$. For a subgraph $H\subseteq G$, let $\varphi|_H$ denote the restriction of $\varphi$ to $H$, i.e., $\varphi|_H(f)=\varphi(f)$ for every $f \in E(H)$. For a vertex $u\in V(G)$, let $\varphi(u)$ denote the set of colors assigned to edges incident with $u$, and $\overline{\varphi}(u)=[k]\backslash\varphi(u)$, that is, the set of colors not assigned to any edge incident with $u$. We call $\varphi(u)$ the set of colors {\bf present} at $u$ and $\overline{\varphi}(u)$ the set of colors {\bf missing} at $u$. Clearly, $|\varphi(u)|+|\overline{\varphi}(u)|=k$ for each vertex $u\in V(G)$. For a vertex set $X\subseteq V(G)$, we define $\overline{\varphi}(X)=\cup_{u\in X} \overline{\varphi}(u)$. The set $X$ is called {\bf elementary} with respect to $\varphi$ or simply {\bf $\varphi$-elementary} if $\overline{\varphi}(u)\cap\overline{\varphi}(v)=\emptyset$ for any two distinct vertices $u,v\in X$. For a color $\alpha\in[k]$, let $E_{\varphi,\alpha}(G)$ denote the set of edges colored with $\alpha$. Let $\alpha,\beta\in[k]$ be two distinct colors, and let $H$ be the subgraph induced by $E_{\varphi,\alpha}(G)$ and $E_{\varphi,\beta}(G)$.
Clearly, every component of $H$ is either a path or an even cycle, each of which is referred to as an $(\alpha,\beta)$-chain of $G$. A \textbf{Kempe change} is the operation of interchanging the colors $\alpha$ and $\beta$ of every edge on an $(\alpha,\beta)$-chain, while keeping the colors of all other edges unchanged. Notably, a Kempe change performed on a whole $(\alpha,\beta)$-chain yields another proper edge-coloring. For a vertex $u$ of $G$, we denote by $P_u(\alpha,\beta,\varphi)$ or simply $P_u(\alpha,\beta)$ the unique $(\alpha,\beta)$-chain containing the vertex $u$. By convention, we allow isolated vertices $u\in V(G)$ as chains with zero edges, and note that a multigraph may contain even cycles of length 2. For any two distinct vertices $u,v\in V(G)$, the two chains $P_u(\alpha,\beta,\varphi)$ and $P_v(\alpha,\beta,\varphi)$ are either identical or disjoint, which is an important fact in our proofs. In particular, we may change the color of one edge in our proofs and adopt the notation $e:\alpha\to\beta$ to mean that we recolor the edge $e$ from $\alpha$ to $\beta$.

\noindent\textbf{Definition 1} (Multi-fan). Let $G$ be a multigraph with an edge  $e_1\in E_G(x,y_1)$, and an edge-coloring $\varphi\in\mathcal{C}^k(G-e_1)$ for some integer $k\geq \Delta(G)$. A multi-fan at $x$ with respect to $e_1$ and $\varphi$ is a sequence $F = (x,e_1,y_1,\dots,e_p,y_p)$ with $p\geq 1$ consisting of edges $e_1,\dots,e_p$ and vertices $x,y_1,\dots,y_p$ satisfying the following two conditions:
	
\begin{itemize}
\item[F1.] The edges $e_1,\dots,e_p$ are distinct, and $e_i\in E_G(x,y_i)$ for $1\leq i\leq p$.
\item[F2.] For every edge $e_i$ with $2\leq i\leq p$, there exists a vertex $y_j$ with $1\leq j\leq i-1$ such that $\varphi(e_i)\in\overline{\varphi}(y_j)$.
\end{itemize}
Here $V(F)$ denotes the set of vertices of $F$. For any $z\in V(F) \setminus\{x\}$, let $\mu_F(x,z)$ denote the number of edges in $F$ connecting $x$ and $z$. We note that the leaves $y_1,\dots,y_p$ of a multi-fan are allowed to be repeated, i.e., $y_i$ and $y_j$ may be identical for $i\neq j$. A multi-fan $F$ with root $x$, initial edge $e_1\in E_G(x,y_1)$, and with respect to a proper edge-coloring $\varphi$ is called maximal if there is no multi-fan properly containing $F$ that shares the same root $x$, the same initial edge $e_1$, and the same coloring $\varphi$.

\noindent\textbf{Definition 2} (Kierstead path). Let $G$ be a multigraph with an edge $e_1\in E_G(y_0,y_1)$, and an edge-coloring $\varphi\in\mathcal{C}^k(G-e_1)$ for some integer $k\geq \Delta(G)$. A Kierstead path with respect to $e_1$ and $\varphi$ is a sequence $K=(y_0,e_1,y_1,\dots,e_p,y_p)$ with $p\geq 1$ consisting of edges $e_1,\dots,e_p$ and vertices $y_0,y_1,\dots,y_p$ satisfying the following two conditions:
	
\begin{itemize}
\item[K1.] The vertices $y_0,\dots,y_p$ are distinct,  and $e_i\in E_G(y_{i-1},y_i)$ for $1\leq i\leq p$.
\item[K2.] For every edge $e_i$ with $2\leq i\leq p$, there exists a vertex $y_j$ with $0\leq j\leq i-2$ such that $\varphi(e_i)\in\overline{\varphi}(y_j)$.
\end{itemize}
	
It is easy to see that a Kierstead path consisting of at most three vertices is also a multi-fan. The two-vertex case is immediate. In the three-vertex case, it can be reordered as $(y_1,e_1,y_0,e_2,y_2)$, a multi-fan rooted at $y_1$.

\begin{lemma}[{[\citealp{Stiebitz2012}, Theorem 2.1]}]\label{L2.1}
Let $G$ be a multigraph with $\chi'(G)=k+1$ for an integer $k\ge\Delta(G)$, and let $e_1\in E_G(x,y_1)$ be a critical edge of $G$. Furthermore, let $F=(x,e_1,y_1,\ldots,e_p,y_p)$ be a multi-fan at $x$ with respect to $e_1$ and $\varphi\in\mathcal{C}^k(G-e_1)$. Then the following statements hold:
		
\begin{enumerate}
\item[\upshape(1)]
$V(F)$ is $\varphi$-elementary;
\item[\upshape(2)]If $\alpha\in\overline{\varphi}(x)$ and $\beta\in\overline{\varphi}(y_i)$ for $1\le i\le p$, then there is an $(\alpha,\beta)$-chain with respect to $\varphi$ having endvertices $x$ and $y_i$;
\item[\upshape(3)]If $F$ is a maximal multi-fan at $x$ with respect to $e_1$ and $\varphi$, then $|V(F)|\ge 3$ and $\sum_{z\in V(F)\setminus\{x\}}(d_G(z)+\mu_F(x,z)-k)=2$.
\end{enumerate}
\end{lemma}

The following lemma concerns the double fan formed by taking the union of two multi-fans at two adjacent vertices. This construction was later extended to the $e$-fan, and analogous results were obtained \cite{CCQ2023}.

\begin{lemma}[{[\citealp{Stiebitz2012}, Theorem 2.17]}]\label{L2.2}
Let $G$ be a multigraph with $\chi'(G)=k+1$ for an integer $k\ge\Delta(G)$, let $e\in E_G(x,y)$ be a critical edge of $G$, and let $\varphi\in\mathcal{C}^k(G-e)$. Furthermore, denote the double fan by $F_x\cup F_y$, where $F_x$ and $F_y$ are the multi-fans at $x$ and $y$, respectively, with respect to $e$ and $\varphi$. If $\min\{d_G(x),d_G(y)\}<k$, then $V(F_x)\cup V(F_y)$ is $\varphi$-elementary.
\end{lemma}

\begin{lemma}[{[\citealp{Stiebitz2012}, Theorem 3.3]}]\label{L2.3}
Let $G$ be a multigraph with $\chi'(G)=\Delta(G)+1$. Let $e_1\in E_G(y_0,y_1)$ be a critical edge, and let $\varphi \in\mathcal{C}^\Delta(G-e_1)$. Let $K=(y_0,e_1,y_1,e_2,y_2,e_3,y_3)$ be a Kierstead path with respect to $e_1$ and $\varphi$. If $\min\{d_G(y_1),d_G(y_2)\}<\Delta(G)$, then $V(K)$ is $\varphi$-elementary.
\end{lemma}
	
\begin{lemma}[Vizing's Adjacency Lemma~\cite{Vizing1965}]\label{L2.4}
Let $G$ be a simple graph with $\chi'(G)=\Delta(G)+1$, and let $e=xy$ be a critical edge of $G$. Then $x$ is adjacent to at least $\Delta(G)-d_G(y)+1$ vertices $z\neq y$ such that $d_G(z)=\Delta(G)$.
\end{lemma}

\begin{lemma}[Hilton–Zhao Core Theorem~\cite{CCJS2024}]\label{L2.5}
Let $G$ be a connected simple graph with $\Delta(G)\ge4$ and $\Delta(G_{\Delta})\le2$. Then $G$ is of Class~2 if and only if $G$ is overfull.
\end{lemma}

\section{Proofs of Theorems \ref{T1} and \ref{T2}}

We begin with the proof of Theorem \ref{T1}.
\begin{proof}	
(1) $\Rightarrow$ (2) :
Since $\chi'_{loc}(G)=\chi'(G)=k$, we can choose a locally bounded proper edge-coloring $\varphi:E(G)\to[k]$.
For every integer $t\in[k]$, if $e\in E(G[D_t])$, then $\varphi(e)\le b_G(e)\le t$. Thus all edges in $E(G[D_t])$ can only use colors from $[t]$. Hence $|\{\varphi(e):e\in E(G[D_t])\}|\le t$.

(2) $\Rightarrow$ (3) : Let $\varphi$ be a proper $k$-edge-coloring satisfying \textup{(2)}, and let $M_i=\varphi^{-1}(i)$ for $i\in[k]$. These color classes form a partition of $E(G)$ into matchings. Since $\chi'(G)=k$, every $M_i$ is nonempty. Hence $b(M_i)=\min_{e\in M_i}b_G(e)$ is well defined, and we have $b(M_i)\le t$ if and only if at least one edge in $M_i$ belongs to $G[D_t]$.
Therefore, the condition $|\{\varphi(e):e\in E(G[D_t])\}|\le t$ is equivalent to $|\{i:b(M_i)\le t\}|\le t$.
Relabel the matchings so that, writing $b_i=b(M_i)$, we have $b_1\le\cdots\le b_k$. Suppose to the contrary that $b_j<j$ for some $j\in[k]$, and put $t=b_j$.
Thus $1\le t<j\le k$. At least the first $j$ color classes satisfy $b(M_i)\le t$, that is, $|\{i:b(M_i)\le t\}|\ge j>t$.
This contradicts the condition $|\{i:b(M_i)\le t\}|\le t$.
Consequently, $b_i\ge i$ for all $i\in[k]$.
		
(3) $\Rightarrow$ (1) : Let the matchings be indexed as in \textup{(3)} with $b_i=b(M_i)$ and $b_1\le\dots\le b_k$. Define $\psi:E(G)\to[k]$ by $\psi(e)=i$ for every edge $e\in M_i$. This is a proper edge-coloring because the $M_i$ form a partition into matchings. For every $e\in M_i$, we have $\psi(e)=i\le b_i\le b_G(e)$. Thus $\psi$ is a locally bounded proper edge-coloring and $k=\chi'(G)\le\chi'_{loc}(G)\le k$. We obtain $\chi'_{loc}(G)=\chi'(G)=k$.	
\end{proof}
	
\medskip
We next turn to the proof of Theorem \ref{T2}.	
		
\begin{proof}	
\medskip
\noindent(1)
Let $G$ be a regular multigraph with degree $r\ge 1$.  Choose a proper $k$-edge-coloring $\varphi$ of $G$ such that $(|M^\varphi_k|,|M^\varphi_{k-1}|,\dots,|M^\varphi_1|)$ is lexicographically minimum, where $M^\varphi_i:=\{e\in E(G):\varphi(e)=i\}$ for $i\in [k]$. Suppose that there exists an edge $e\in E_G(u,v)$ with color $t\in [k]$ that does not satisfy the local bound. By the regularity, $t>r+\max\{\mu_G(u),\mu_G(v)\}$, and so $t-1\ge r+\mu_G(u)>r$. Let $H$ be the component containing $e$ in the subgraph induced by the edge set $\{f:\varphi(f)\in[t-1]\}\cup\{e\}$. Suppose that $H$ admits a proper $(t-1)$-edge-coloring. Then recoloring edges in $E(H)$ using colors from $[t-1]$ can yield a new proper $k$-edge-coloring $\varphi'$ of $G$ since all edges in $E(G)\setminus E(H)$ incident with a vertex of $H$ have colors at least $t$. Now we have $|M^{\varphi'}_t|=|M^\varphi_t|-1$ and $|M^{\varphi'}_j|=|M^\varphi_j|$ for $j>t$ if such $j$ exists, contradicting the choice of $\varphi$. Hence $\chi'(H)=t$. Since $\varphi|_{H-e}$ is a proper $(t-1)$-edge-coloring of $H-e$, $e$ is a critical edge of $H$. Note that $\Delta(H)\le r<t-1$. Let $F$ be a maximal multi-fan at $u$ in $H$ with respect to $e$ and $\varphi|_{H-e}$. By Lemma~\ref{L2.1}(3), we obtain that $\sum_{z\in V(F)\setminus\{u\}}\big[d_H(z)+\mu_H(u,z)-(t-1)\big]\ge \sum_{z\in V(F)\setminus\{u\}}\big[d_H(z)+\mu_F(u,z)-(t-1)\big]=2$. However, each term in the sum satisfies $d_H(z)+\mu_H(u,z)-(t-1)\le r+\mu_G(u)-(t-1)\le 0$, which yields a contradiction.
		
Therefore, all edges in $E(G)$ satisfy the local bound, and so the chosen proper $k$-edge-coloring $\varphi$ is locally bounded. Consequently, we have $k=\chi'(G)\le\chi'_{loc}(G)\le k$, which implies $\chi'_{loc}(G)=\chi'(G)$ as desired.

\medskip
Recall that $B^\varphi_i$ denotes the set of non-local $i$-edges with $i\in[k]$.

\noindent(2)\label{2}
Choose a proper $k$-edge-coloring $\varphi$ of $G$ such that $(|B^\varphi_k|,|B^\varphi_{k-1}|,\ldots,|B^\varphi_1|)$ is lexicographically minimum. Notice that $|B^\varphi_1|=|B^\varphi_2|=0$ since every edge satisfies $b_G(e) \ge2$. We first show that $|B^\varphi_4|=|B^\varphi_3|=0$. Suppose that $e\in E_G(u,v)$ is a non-local 4-edge with $\varphi(e)=4$. It follows that $d_G(u),d_G(v)\le2$ since $4>\max\{d_G(u)+\mu_G(u), d_G(v)+\mu_G(v)\}$. Thus, we obtain that $|(\varphi(u)\cup\varphi(v))\cap[3]|\le2$ and $(\overline{\varphi}(u)\cap\overline{\varphi}(v))\cap[3]\neq\emptyset$. Recoloring $e$ with a common missing color of $u$ and $v$  decreases $|B_4^\varphi|$, contradicting the choice of $\varphi$. Suppose that $e\in E_G(u,v)$ is a non-local 3-edge with $\varphi(e)=3$. It follows that $d_G(u)=d_G(v)=1$ since $3>\max\{d_G(u)+\mu_G(u), d_G(v)+\mu_G(v)\}$, and so $e$ has no adjacent edges. Thus recoloring $e$ with color $1$ decreases $|B_3^\varphi|$, also contradicting the choice of $\varphi$. Thus, $|B^\varphi_4|=|B^\varphi_3|=0$ as desired.

Now suppose that there exists a non-local $t$-edge $e\in E_G(u,v)$ with $t\ge5$. Since $t>b_G(e) =\max\{d_G(u)+\mu_G(u),d_G(v)+\mu_G(v)\}$, we have $d_G(u)<t-1$ and $d_G(v)<t-1$. Let $H$ be the component containing $e$ in the subgraph induced by the edge set $\{f:\varphi(f)\in [t-1]\}\cup\{e\}$. We have $\Delta(H)\le t-1$ since $d_G(u),d_G(v)<t-1$ and other vertices in $V(H)$ are incident with at most $t-1$ edges. Suppose that $H$ admits a proper $(t-1)$-edge-coloring. Then recoloring edges in $E(H)$ yields a new proper $k$-edge-coloring $\varphi'$ of $G$ since all edges in $E(G)\setminus E(H)$ incident with a vertex of $H$ have colors at least $t$. Now $|B^{\varphi'}_t|=|B^{\varphi}_t|-1$ and $|B^{\varphi'}_j|=|B^{\varphi}_j|$ for every $j\in [k]$ with $j>t$, contradicting the choice of $\varphi$. Hence $\chi'(H)=t$, and $e$ is a critical edge of $H$ since $\varphi|_{H-e}$ is a proper $(t-1)$-edge-coloring of $H-e$.

Suppose first that $e\in E_G(u,v)$  is a non-local $t$-edge with $t\ge6$. If $d_G(u),d_G(v)\le3$, then $d_{G-e}(u)+d_{G-e}(v)\le4$. Hence $u$ and $v$ have at least one common missing color in $[5]$, say $t'$. We can recolor $e$ with the color $t'$ to decrease $|B^\varphi_t|$, contradicting the choice of $\varphi$. Therefore, assume without loss of generality that $d_G(u)=4$. Note that the condition for a non-local $t$-edge gives $4+\mu_G(u)\le t-1$. Let $F$ be a maximal multi-fan at $u$ in $H$ with respect to $e$ and $\varphi|_{H-e}$. By Lemma~\ref{L2.1}(3), we obtain that $\sum_{z\in V(F)\setminus\{u\}}[d_H(z)+\mu_H(u,z)-(t-1)]\ge\sum_{z\in V(F)\setminus\{u\}}[d_H(z)+\mu_F(u,z)-(t-1)]=2$. However, each term in the sum satisfies $d_H(z)+\mu_H(u,z)-(t-1)\le4+\mu_G(u)-(t-1)\le0$, a contradiction.

Finally, we consider the case $t=5$. Suppose that $e\in E_G(u,v)$ is a non-local 5-edge. Then we have $d_G(u)+\mu_G(u)\le4$ and $d_G(v)+\mu_G(v)\le4$. Thus $d_G(u),d_G(v)\le 3$. If $\min\{d_G(u), d_G(v)\}\le 2$, then $d_{G-e}(u)+d_{G-e}(v)\le 3$. Hence $u$ and $v$ have at least one common missing color in $[4]$. We can recolor $e$ with one common missing color of $u$ and $v$ to decrease $|B_t^\varphi|$, contradicting the choice of $\varphi$. Therefore, we have $d_G(u)=d_G(v)=3$ and $\mu_G(u)=\mu_G(v)=1$. We claim that $\varphi(u)$ and $\varphi(v)$ together occupy all four distinct colors in $[4]$. Otherwise, we can recolor $e$ with a common missing color of $u$ and $v$ in $[4]$ to decrease $|B_t^\varphi|$, contradicting the choice of $\varphi$. Without loss of generality, we may assume that $\varphi(u)=\{1,2,5\}$ and $\varphi(v)=\{3,4,5\}$. Let $\varphi(e_1)=\alpha\in\{1,2\}$, where $e_1\in E_G(u,w)$. Note that in the subgraph $H$, $(u,e,v,e_1,w)$ forms a multi-fan at $u$ with respect to $e$ and  $\varphi|_{H-e}$. By Lemma~\ref{L2.1}(1), we have $(\overline{\varphi|_{H-e}}(u)\cup\overline{\varphi|_{H-e}}(v))\cap\overline{\varphi|_{H-e}}(w)=\emptyset$, yielding $d_H(w)=\Delta(H)=4$. Since $\Delta(G)\le4$, we have $d_G(w)=4$, $\varphi(w)=\{1,2,3,4\}$, and so $5\in \overline{\varphi}(w)$. Now simultaneously recolor $e:5\to\alpha$ and $e_1:\alpha\to5$. We obtain a new proper edge-coloring $\varphi''$ of $G$, and $\varphi''(e_1)=5\le d_G(w)+\mu_G(w)\le b_G(e_1)$. Note that $e$ and $e_1$ are not non-local edges under $\varphi''$. Thus $|B^{\varphi''}_5|=|B^\varphi_5|-1$, while $|B^{\varphi''}_j|=|B^\varphi_j|$ ($j\ge6$) remains unchanged, a contradiction.

Therefore, the chosen proper $k$-edge-coloring $\varphi$ is locally bounded, and so $k=\chi'(G)\le\chi'_{loc}(G)\le k$, which implies $\chi'_{loc}(G)=\chi'(G)$ as desired.

\medskip
\noindent(3)
Suppose that the conclusion fails, so $\chi'_{loc}(G)>\chi'(G)$.
Choose a proper $k$-edge-coloring $\varphi$ of $G$ that minimizes the total excess.
That is, $\varphi$ satisfies
\[
\Theta(\varphi):=\sum_{e\in E(G)}(\varphi(e)-b_G(e))^+=\min_{\varphi'\in \mathcal{C}^k(G)}\Theta(\varphi'),
\]
where $(\varphi(e)-b_G(e))^+=\max\{\varphi(e)-b_G(e),0\}$.
Since $\chi'_{loc}(G)>\chi'(G)=k$, there exists at least one non-local edge $e\in E_G(u,v)$ with $\varphi(e)=t>b_G(e)$, where $t\in [k]$. Hence $\Theta(\varphi)>0$ and $d_G(u)+d_G(v)\le b_G(e)+\mu_G(u,v)$. Since $\varphi(e)=t$ and $t\notin[t-1]$, only the remaining $d_G(u)-1$ edges incident to $u$ and $d_G(v)-1$ edges incident to $v$ can possibly occupy colors in $[t-1]$. Since every color assigned to an edge between $u$ and $v$ appears at both endvertices, we have $|(\varphi(u)\cup\varphi(v))\cap[t-1]|\le d_G(u)+d_G(v)-\mu_G(u,v)-1\le b_G(e)-1\le t-2$. Thus vertices $u$ and $v$ have at least one common missing color in $[t-1]$. Choose a common missing color $t'\in[t-1]$. We recolor the edge $e: t\rightarrow t'$ and keep the colors of all other edges unchanged. Let the new proper $k$-edge-coloring be denoted by $\varphi'$. Then we have $(\varphi'(e)-b_G(e))^+<(\varphi(e)-b_G(e))^+$, and so $\Theta(\varphi')\leq\Theta(\varphi)-1$, contradicting the choice of $\varphi$. Hence $\chi'_{loc}(G)=\chi'(G)$ as desired.
\end{proof}

\section{Proof of Theorem \ref{T3}}

\begin{proof}
If $G$ is of Class $2$, then by equation~\eqref{eq:basic-bounds} we obtain that $\Delta(G)+1=\chi'(G)\le\chi'_{loc}(G)\le\Delta(G)+1$. Thus $\chi'_{loc}(G)=\chi'(G)$ as desired.
Now assume that $G$ is of Class $1$, i.e., $\chi'(G)=\Delta(G)$. The case $\Delta(G)\leq 4$ has been proved by Theorem~\ref{T2}(2). Hence, we may assume without loss of generality that $\Delta(G)\ge 5$.
		
Let $C_a \vee \overline{K}_b$ denote the join of the cycle $C_a$ with $a$ vertices and the edgeless graph $\overline{K}_b$ with $b$ vertices, which is constructed from the disjoint union of $C_a$ and $\overline{K}_b$ by adding all edges between the vertex set of $C_a$ and the vertex set of $\overline{K}_b$.
		
\begin{Claim}
\label{claim 1}
$G$ contains no $r$-critical subgraph with $r\ge6$, and every $5$-critical subgraph of $G$ is isomorphic to $C_4\vee\overline{K}_3$.
\end{Claim}
		
\begin{proof}
Let $H\subseteq G$ be an $r$-critical subgraph with  $\chi'(H)=\Delta(H)+1$ and $\Delta(H)=r\geq 5$.  Notice that $H$ is connected by criticality, and the core of $H$ is $H_r$ with $V(H_r)=\{x\in V(H):d_H(x)=r\}$. We have $\Delta(H_r)\leq 2$ since $V(H_r)\subseteq V_{\ge5}(G)$ and $\Delta(G[V_{\ge5}(G)])\leq 2$. It follows that $H$ is overfull by Lemma~\ref{L2.5}. By Lemma~\ref{L2.4}, every $r$-vertex has at least two $r$-neighbors in $H$. Therefore, $H_r$ is regular with degree 2. For any $y\in V(H)\setminus V(H_r)$, we can choose an $r$-neighbor $z\in V(H_r)$ by Lemma~\ref{L2.4}. Note that $z$ has at least $r-d_H(y)+1$ $r$-neighbors by Lemma~\ref{L2.4}. Since $d_H(y)<r$, we have $2\le r-d_H(y)+1\le |N_{H_r}(z)|=2$, hence $d_H(y)=r-1$.
			
If $r\geq 6$, then all vertices in $H$ have degree at least 5 in $G$. Hence $H\subseteq G[V_{\ge5}(G)]$. Since $\Delta(G[V_{\ge5}(G)])\leq 2$, we have $\Delta(H)\leq 2$, which contradicts $\Delta(H)=r\geq 6$. Thus there is no $r$-critical subgraph with maximum degree $r\geq 6$ in $G$.
			
Now assume $r=5$. The degrees of vertices in $H$ are only 4 or 5. Write $n=|V(H)|$. Since $H$ is overfull, we have $|E(H)|>5\lfloor\frac{n}{2}\rfloor$. For every $e\in E(H)$, there exists a proper $5$-edge-coloring of $H-e$ since $H$ is $5$-critical. Thus $|E(H)|-1\le 5\lfloor\frac{n}{2}\rfloor$, and so $|E(H)|=5\lfloor\frac{n}{2}\rfloor+1$. If $n$ is even, then we have $2|E(H)|=5n+2$, contradicting $2|E(H)|\le 5n$. Hence $n$ is odd, and we obtain $|E(H)|=\frac{5(n-1)}{2}+1$. Let $n_4$ and $n_5$ be the number of 4-vertices and the number of 5-vertices in $H$, respectively. We have $n_4+n_5=n$ and $4n_4+5n_5=2|E(H)|=5(n-1)+2$. It follows that $n_4=3$ and $n_5=n-3$.
			
We next characterize the structure of $H$. We have $n_5\geq 3$ since the core $H_5$ is regular with degree 2. Each 5-vertex has exactly two 5-neighbors, and its remaining three neighbors are precisely the three 4-vertices. Consequently, every 4-vertex must be adjacent to all 5-vertices, and so $n_5\leq 4$. Now we obtain $3\le n_5\le4$. It follows that $n_5=4$ since $n$ is odd and $n_5=n-3$.
Therefore, we have the core $H_5\cong C_4$, and the three 4-vertices form an independent set $\overline{K}_3$, i.e., $H\cong C_4 \vee \overline{K_3}$.
\end{proof}

Note that $B^\varphi_i:=\{xy:\varphi(xy)=i,\max\{d_G(x),d_G(y)\}\leq i-2\}$ with $i\in[\Delta]$.

\begin{Claim}
\label{claim 2}
We may assume $|B^\varphi_t|=0$ for every integer $t$ with $6\le t\le\Delta$.
\end{Claim}

\begin{proof}
Assume $\Delta\ge 6$. Choose a proper $\Delta$-edge-coloring $\varphi$ of $G$ that lexicographically minimizes
\begin{equation}
(|B^\varphi_\Delta|,|B^\varphi_{\Delta-1}|,\dots,|B^\varphi_6|).
\label{eq:lexicographically}
\end{equation}
Suppose to the contrary that there exists a non-local $t$-edge $xy$ with $t\geq 6$. We have $t>\max\{d_G(x)+\mu_G(x), d_G(y)+\mu_G(y)\}$, so it follows  that
\begin{equation}
d_G(x),d_G(y)\leq t-2.
\label{eq:degree}
\end{equation}
		
Let $H_1$ be the component containing $xy$ of the induced subgraph with the edge set $\{e:\varphi(e)\in[t-1]\}\cup\{xy\}$.
We have $\Delta(H_1)\le t-1$ since every vertex in $V(H_1)$ is incident with at most $t-1$ edges in $H_1$.
We claim that $\chi'(H_1)=t$, $\Delta(H_1)=t-1$, and $xy$ is a critical edge of $H_1$. Suppose that $H_1$ admits a proper $(t-1)$-edge-coloring. Then recoloring edges in $E(H_1)$ yields a new proper $\Delta$-edge-coloring $\varphi'$ of $G$ since all edges in $E(G)\setminus E(H_1)$ incident with a vertex of $H_1$ have colors at least $t$. Now we have $|B^{\varphi'}_t|=|B^{\varphi}_t|-1$ and $|B^{\varphi'}_j|=|B^{\varphi}_j|$ for $j>t$ if such $j$ exists, contradicting the choice of $\varphi$. Hence $\chi'(H_1)=t$, and $xy$ is a critical edge of $H_1$ since $\varphi|_{H_1-xy}$ is a proper $(t-1)$-edge-coloring of $H_1-xy$. Note that $H_1$ is a simple graph and $\Delta(H_1)\le t-1$.  We get $\Delta(H_1)=t-1$ by Vizing's theorem.
		
Choose a $(t-1)$-critical subgraph $H_2\subseteq H_1$ with $\chi'(H_2)=\chi'(H_1)=t$. We obtain $\Delta(H_2)=\Delta(H_1)=t-1$ by Vizing's theorem. Since $xy$ is a critical edge of $H_1$, $xy$ must be in $E(H_2)$. Otherwise, if $xy\notin E(H_2)$, then $H_2\subseteq H_1-xy$, and so $t=\chi'(H_2)\leq\chi'(H_1-xy)=t-1$, a contradiction.
If $t\ge 7$, then we have $\Delta(H_2)=t-1\ge6$, which contradicts Claim~\ref{claim 1}. If $t=6$, then we have $\Delta(H_2)=t-1=5$. By Claim~\ref{claim 1}, we have that $H_2\cong C_4\vee\overline{K}_3$ and $d_{H_2}(x),d_{H_2}(y)\ge 4$.
Together with~\eqref{eq:degree}, this gives $d_{H_2}(x)=d_{H_2}(y)=4$. However, all 4-vertices are pairwise non-adjacent in $H_2$, which contradicts the fact $xy\in E(H_2)$.
\end{proof}

Let $S_5^\varphi:=\{z\in V(G):d_G(z)=4 \text{ and }5\in{\varphi}(z)\}$.

\begin{Claim}
\label{claim 3}
 We may continue to assume $|B^\varphi_5|=0$.
\end{Claim}

\begin{proof}   	
By Claim \ref{claim 2},
we may continue to choose a proper $\Delta$-edge-coloring $\varphi$ of $G$ such that $|B^\varphi_\Delta|=|B^\varphi_{\Delta-1}|=\cdots=|B^\varphi_6|=0$ and $(|B_5^\varphi|,-|S_5^\varphi|)$ is lexicographically minimum.
		
Suppose to the contrary that there exists a non-local $5$-edge $xy$, so that $|B_5^\varphi|>0$. Then $d_G(x)\le 3$ and $d_G(y)\le 3$. If $\min\{d_G(x), d_G(y)\}<3$, then $x$ and $y$ have a common missing color in $[4]$, say $\alpha$. We may directly recolor $xy$ with the color $\alpha$, so that $|B_5^\varphi|$ is decreased by one under the new edge-coloring, which contradicts the choice of $\varphi$. Therefore, $d_G(x)=d_G(y)=3$. We claim that $\varphi(x)$ and $\varphi(y)$ together occupy all four distinct colors in $[4]$. Otherwise, we can recolor $xy$ with a common missing color of $x$ and $y$ in $[4]$ to decrease $|B_5^\varphi|$, contradicting the choice of $\varphi$. Thus, $\varphi(x)\cup\varphi(y)\supseteq[4]$ as desired. Without loss of generality, we may assume that $\varphi(xy)=5$, $\varphi(x)=\{1,2,5\}$, and $\varphi(y)=\{3,4,5\}$.
		
Let $H$ be the component containing $xy$ of the induced subgraph with edge set $\{e\in E(G):\varphi(e)\in [4]\}\cup\{xy\}$. We claim that $\chi'(H)=5$, $\Delta(H) = 4$, and $xy$ is a critical edge of $H$. Suppose that $H$ admits a proper $4$-edge-coloring. Then recoloring edges in $E(H)$ yields a new proper $\Delta$-edge-coloring $\varphi'$ of $G$ since all edges in $E(G)\setminus E(H)$ incident with a vertex of $H$ have colors at least $5$. Now $|B^{\varphi'}_5|=|B^{\varphi}_5|-1$, contradicting the choice of $\varphi$. Hence $\chi'(H)=5$, and $xy$ is a critical edge of $H$ since $\varphi|_{H-xy}$ is a proper $4$-edge-coloring of $H-xy$. Note that $H$ is a simple graph and $\Delta(H)\le 4$.  We get $\Delta(H)=4$ by Vizing's theorem.

Let $\varphi_1:=\varphi|_{H-xy}$. Obviously, $d_H(x)=d_H(y)=3$, $\overline{\varphi}_1(x)=\{3,4\}$, and $\overline{\varphi}_1(y)=\{1,2\}$. Let $F_x\cup F_y$ be the double fan with respect to $xy$ and $\varphi_1$. Since $d_H(x)=d_H(y)=3<\Delta(H)=4$, it follows that $V(F_x)\cup V(F_y)$ is $\varphi_1$-elementary by Lemma \ref{L2.2}.

\begin{Subclaim}
\label{subclaim31}
Every vertex in $V(H)\setminus\{x,y\}$ at distance $1$ or $2$ from the set $\{x,y\}$ in $H$ has degree $4$ in $H$.

\end{Subclaim}
\begin{proof}
Since $\varphi_1(x)=\overline{\varphi}_1(y)$ and $\varphi_1(y)=\overline{\varphi}_1(x)$, we have $N_H(x)\cup N_H(y)=V(F_x)\cup V(F_y)$. For any $z\in(N_H(x)\cup N_H(y))\setminus\{x,y\}$, we have $d_H(z)=\Delta(H)=4$ since $\overline{\varphi}_1(x)\cup\overline{\varphi}_1(y)=[4]$ and $V(F_x)\cup V(F_y)$ is $\varphi_1$-elementary.
		
Let $z\in N_H(x)\setminus\{y\}$ and $u\in N_H(z)\setminus\{x,y\}$. We obtain a Kierstead path $K:=(y,yx,x,xz,z,zu,u)$, since $\varphi_1(xz)\in\overline{\varphi}_1(y)\text{ and } \varphi_1(zu)\in\overline{\varphi}_1(x)\cup\overline{\varphi}_1(y)=[4]$.
Since $d_H(x)=3<\Delta(H)=4$, $V(K)$ is $\varphi_1$-elementary by Lemma~\ref{L2.3}. It follows that $d_H(u)=\Delta(H)=4$ since $\overline{\varphi}_1(x)\cup\overline{\varphi}_1(y)=[4]$.
The same argument applies to $z\in N_H(y)\setminus\{x\}$ and $u\in N_H(z)\setminus\{x,y\}$.
\end{proof}		
	
\begin{Subclaim} 	
\label{subclaim32}
$d_G(z)\ge 5$ and $[5]\subseteq\varphi(z)$ for every $z\in (N_H(x)\cup N_H(y))\setminus\{x,y\}$.
\end{Subclaim}

\begin{proof}
Let $z\in N_H(x)\setminus\{y\}$ and $\varphi(xz)=\alpha\in\{1,2\}$. By Subclaim~\ref{subclaim31}, $d_G(z)\ge d_H(z)=4$ and $\varphi_1(z)=[4]$. Suppose  that $5\in\overline{\varphi}(z)$. Since $\alpha\in\overline{\varphi}(y)$, we can simultaneously recolor $xy:5\to\alpha$ and $xz:\alpha\to5$, which yields a new proper $\Delta$-edge-coloring  of $G$ without forming  new non-local edges and forces $|B_5^\varphi|$ to decrease by one, contradicting the choice of $\varphi$. Thus $[5]\subseteq\varphi(z)$ and $d_G(z)\ge 5$. The same argument applies to vertices in $N_H(y)\setminus\{x\}$.
\end{proof}

\begin{Subclaim}
\label{subclaim33}
For $\alpha\in\{3,4\}$ and $\beta\in\{1,2\}$, $P_x(\alpha,\beta,\varphi)=P_y(\alpha,\beta,\varphi)$.
\end{Subclaim}

\begin{proof}
Note that $(x,xy,y)$ is a multi-fan with respect to $xy$ and $\varphi_1$ in $H$. For $\alpha\in\overline{\varphi}_1(x)=\{3,4\}$ and $\beta\in\overline{\varphi}_1(y)=\{1,2\}$,  there exists an $(\alpha,\beta)$-chain  with endvertices $x$ and $y$ with respect to $\varphi_1$ by Lemma~\ref{L2.1}(2). Since $\varphi(x)=\{1,2,5\}$ and $\varphi(y)=\{3,4,5\}$, we have $\alpha\in\overline{\varphi}(x)$ and $\beta\in\overline{\varphi}(y)$ under $\varphi$. Thus the endvertices of the $(\alpha,\beta)$-chain are also $x$ and $y$ in $G$ under $\varphi$. Therefore, $P_x(\alpha,\beta,\varphi)=P_y(\alpha,\beta,\varphi)$.
\end{proof}
		
\begin{Subclaim}
\label{subclaim34}
Assume that $z\in (N_G(x)\cup N_G(y))\setminus\{x,y\}$ and $u\in N_G(z)\setminus\{x,y\}$. Neither of the following configurations can
occur:
\[
\varphi(xz)=\alpha\in\{1,2\},\varphi(zu)=\beta\in\{3,4\},d_G(u)=4, \text{and }\varphi(u)=[4],
\]
or
\[
\varphi(yz)=\alpha\in\{3,4\},\varphi(zu)=\beta\in\{1,2\},d_G(u)=4, \text{and }\varphi(u)=[4].
\]
\end{Subclaim}	
		
\begin{proof}
Suppose to the contrary that the first configuration is satisfied. Let $v$ be an arbitrary vertex satisfying $d_G(v)=4$ and $\varphi(v)=[4]$, so $5\in\overline{\varphi}(v)$. Let $\gamma\in [4]$ and define $P:=P_v(\gamma,5,\varphi)$, where $v$ is an endvertex of $P$. Denote the other endvertex of $P$ by $v'$. Let $e_1=vv_1$ and $e_2=v'v_2$, where $v_1,v_2\in V(P)$. Note that $\varphi(e_1)=\gamma$. When $|E(P)|=1$, we have $e_1=e_2$. We claim that $P$ contains no non-local 5-edge. If $|E(P)|=1$, then the claim holds trivially since $\varphi(e_1)=\gamma$. Assume now that $|E(P)|\ge 2$. If $P$ contains a non-local 5-edge $e_3=ww'$, then it follows that $d_G(w)\le 3$ and $d_G(w')\le 3$ since $5>\max\{d_G(w)+\mu_G(w), d_G(w')+\mu_G(w')\}$. If $e_3$ is an internal edge of $P$, then we have $\gamma\in\varphi(w)\cap\varphi(w')$. Excluding this non-local 5-edge $e_3$, we have $|(\varphi(w)\cup\varphi(w'))\cap [4]|\le 2+2-1=3$, which implies $|\overline{\varphi}(w)\cap \overline{\varphi}(w')\cap [4]|\ge 1$. Thus, we can recolor $e_3$ with a common missing color of $w$ and $w'$, so that $|B_5^\varphi|$ is decreased by one under the new edge-coloring, which contradicts the choice of $\varphi$. Therefore, there are no non-local 5-edges in the interior of $P$. If $e_3$ is incident to an endvertex, then $e_3=e_2$ since $\varphi(e_1)=\gamma$. Apply a Kempe change on $P$, and let $P'$ denote the path obtained after swapping colors $\gamma$ and 5 along $P$. After the interchange, the internal 4-vertices still see color 5 and $e_1$ is colored with 5. Since $d_G(v)=4$, $e_1$ is not a non-local 5-edge. Meanwhile, $e_2$ is no longer colored with 5, reducing the number of non-local 5-edges by one. If new non-local 5-edges appear in the interior of $P'$, then both endvertices of every non-local 5-edge have degree at most $3$. Moreover, the color $\gamma$ is present at both endvertices. Therefore, the two endvertices together occupy at most three colors from $[4]$, which implies the existence of a common missing color. Hence we can recolor these new non-local 5-edges with a common missing color in $[4]$. Since these new non-local 5-edges form a matching, these recolorings do not interfere with one another and do not change the edges of color $5$ at the 4-vertices. Therefore, after the interchange, $|B_5^\varphi|$ decreases, contradicting the choice of $\varphi$. We conclude that $P$ contains no non-local 5-edge. Thus, we obtain $xy\notin E(P)$.

We apply a Kempe change on $P$ and recolor the newly created internal non-local 5-edges to obtain a new coloring $\varphi'$. Hence $\gamma\in\overline{\varphi'}(v)$ and $5\in\varphi'(v)$. The edge $e_1$ is not a non-local 5-edge since $d_G(v)=4$. Therefore, $e_2$ is the only possible newly created non-local 5-edge after recoloring. So
\begin{equation}
|B^{\varphi'}_5|=|B^\varphi_5|+\varepsilon,\text{ where }\varepsilon\in\{0,1\}.
\label{eq:B}
\end{equation}
For the internal vertices of $P$, the recoloring operations are performed only on non-local 5-edges. Hence the edges colored with 5 that are incident to 4-vertices will not be recolored. When $\varepsilon=1$, we have $d_G(v')\le 3$. Thus, $v'\notin S_5^\varphi$ and $v'\notin S_5^{\varphi'}$. Since $d_G(v)=4$, $\varphi(e_1)=\gamma$, and $\varphi'(e_1)=5$, it follows that $v \notin S_5^\varphi$ and $v\in S_5^{\varphi'}$. Therefore,
\begin{equation}
|S^{\varphi'}_5|=|S_5^\varphi|+1.
\label{eq:S}
 \end{equation}
When $\varepsilon=0$, $|S^{\varphi'}_5|$ does not necessarily increase by one.

We now apply the preceding construction to the original coloring $\varphi$, with $v=u$ and $\gamma=\alpha$, and write $P:=P_u(\alpha,5,\varphi)$. Let $\varphi'$ be the coloring obtained from this particular $P$ by the Kempe change followed by the internal recolorings described above. Since $\varphi(xz)=\alpha$ and $\varphi(xy)=5$, the vertices $x,y,z$ lie in the same
$(\alpha,5)$-chain. It follows that $x,y,z\notin V(P)$ since $xy\notin E(P)$. Hence, the colors of $xy$ and $xz$ remain unchanged under the coloring $\varphi'$. We can simultaneously recolor $xy:5\to\alpha$, $xz:\alpha\to \beta$, and $zu:\beta\to\alpha$ to obtain a new coloring $\varphi''$. If $P$ has an even number of edges, then $\varphi(e_2)=5$ and $\varphi''(e_2)=\alpha$. We have $e_2\notin B_5^\varphi$ and $e_2\notin B_5^{\varphi''}$, whereas $xy\in B_5^\varphi$ but $xy\notin B_5^{\varphi''}$. Therefore, $|B_5^{\varphi''}|=|B_5^\varphi|-1$, contradicting the choice of $\varphi$. If $P$ has an odd number of edges, then $\varphi(e_2)=\alpha$ and $\varphi''(e_2)=5$. If $e_2\notin B_5^{\varphi''}$, then $|B_5^{\varphi''}|=|B_5^\varphi|-1$, contradicting the choice of $\varphi$. If $e_2\in B_5^{\varphi''}$, then $\varepsilon=1$. We obtain $|B_5^{\varphi''}|=|B^{\varphi'}_5|-1=|B_5^\varphi|$. By equation~\eqref{eq:S}, when $\varepsilon=1$, we have $|S_5^{\varphi''}|=|S_5^{\varphi'}|=|S_5^\varphi|+1$, contradicting the choice of $\varphi$.

The second configuration is excluded by interchanging the roles of $x$ and $y$.
\end{proof}
	
Recall that $\varphi(xy)=5$, $\varphi(x)=\{1,2,5\}$, and $\varphi(y)=\{3,4,5\}$. We distinguish two cases according to
whether $N_G(x)\cap N_G(y)$ is empty.
		
\noindent\textbf{Case 1.}
$N_G(x)\cap N_G(y)=\emptyset$.
		
Without loss of generality, assume that $z\in N_G(x)\setminus\{y\}$ and $\varphi(xz)=1$.
We have $[5]\subseteq\varphi(z)$ by Subclaim \ref{subclaim32}.
Now we may define the edges incident to $z$, which are colored with  $2,3,4$, and $5$, as:
\[
\varphi(zu)=2,\varphi(zv)=3,\varphi(zw)=4,\text{ and }\varphi(zp)=5.
\]

By Subclaim~\ref{subclaim31}, we have $d_H(u)=d_H(v)=d_H(w)=4$. It follows that $v\in V_{\ge5}(G)$ and $w\in V_{\ge5}(G)$ by Subclaim~\ref{subclaim34}. Since $z\in V_{\ge5}(G)$ and $\Delta(G[V_{\ge5}(G)])\leq 2$, we have $u,p\notin V_{\ge5}(G)$, i.e., $d_G(u),d_G(p)\leq 4$. Together with $d_H(u)=4$, we obtain that $d_G(u)=4$ and $\varphi(u)=[4]$. Furthermore, we have the following subclaim.
		
\begin{Subclaim}
\label{subclaim35}	
$P_u(2,5,\varphi)=(u,uz,z,zp,p)$, $d_G(p)=4$, and $\varphi(p)=\{1,3,4,5\}$.
\end{Subclaim}
	
\begin{proof}
As $\varphi(u)=[4]$, $5\in\overline{\varphi}(u)$. Note that the path $P_u(2,5,\varphi)$ passes through $z$ and $p$ since $\varphi(zu)=2$ and $\varphi(zp)=5$. From the proof of Subclaim~\ref{subclaim34}, it follows that $xy$ is not on the path $P_u(2,5,\varphi)$ since $xy$ is a non-local 5-edge. By applying a Kempe change on $P_u(2,5,\varphi)$ and recoloring the newly formed internal non-local 5-edges, we obtain a coloring $\varphi'$. By equations~\eqref{eq:B} and~\eqref{eq:S}, $|B^{\varphi'}_5|=|B^\varphi_5|+\varepsilon$ with $\varepsilon\in\{0,1\}$. Furthermore, if $\varepsilon=1$, then $|S^{\varphi'}_5|=|S^\varphi_5|+1$.
		
Under $\varphi'$, let $H'$ be the component containing $xy$ of the induced subgraph with the edge set $\{e\in E(G):\varphi'(e)\in [4]\}\cup\{xy\}$. We claim that $\chi'(H')=5$, $\Delta(H')=4$, and  $xy$ is also a critical edge of $H'$.
Suppose that $H'$ admits a proper $4$-edge-coloring. Then recoloring edges in $E(H')$ yields a new proper $\Delta$-edge-coloring $\varphi''$ of $G$ since all edges in $E(G)\setminus E(H')$ incident with a vertex of $H'$ have colors at least $5$. Now $|B^{\varphi''}_5|=|B^{\varphi'}_5|-1=|B^\varphi_5|+\varepsilon-1$. Since $d_G(x)=d_G(y)=3$, we have $|S^{\varphi''}_5|=|S^{\varphi'}_5|$. When $\varepsilon=0$, we have $|B^{\varphi''}_5|=|B^\varphi_5|-1$; when $\varepsilon=1$, we obtain that $|B^{\varphi''}_5|=|B^\varphi_5|$ and $|S^{\varphi''}_5|=|S^{\varphi'}_5|=|S^\varphi_5|+1$. Both cases contradict the choice of $\varphi$. Therefore, $\chi'(H')=5$
and $xy$ is a critical edge of $H'$ since $\varphi'|_{H'-xy}$ is a proper $4$-edge-coloring of $H'-xy$. Note that $H'$ is a simple graph and $\Delta(H')\le 4$.  We also get $\Delta(H') = 4$ by Vizing's theorem.

After interchanging and recoloring, we obtain that $\varphi'(xz)=1$, $\varphi'(zp)=2$, and $\varphi'(zu)=5$. The path $(y,yx,x,xz,z,zp,p)$ satisfies the conditions for a Kierstead path, and $d_{H'}(x)=3<4$. By Lemma~\ref{L2.3}, the set of vertices $\{y,x,z,p\}$ is $\varphi'|_{H'-xy}$-elementary. Since $3,4\in \overline{\varphi'}(x)$ and $1,2\in\overline{\varphi'}(y)$, we have $\overline{\varphi'}(p)\cap[4]=\emptyset$. We conclude that $d_{H'}(p)=4$. Since $d_G(p)\le 4$, we obtain that $d_G(p)=4$ and $\varphi'(p)=\{1,2,3,4\}$. Therefore, these recoloring operations for non-local 5-edges leave the colors of edges incident to $p$ unchanged. Hence, we have $\varphi(p)=\{1,3,4,5\}$ before interchanging. Consequently, $2\in \overline{\varphi}(p)$, which means that $p$ is an endvertex of the $(2,5)$-chain. Since  $P_u(2,5,\varphi)$ passes through $z$ and $p$, we obtain $P_u(2,5,\varphi)=(u,uz,z,zp,p) $.
\end{proof}		
		
By Subclaim~\ref{subclaim33}, we have the path $P:=P_x(2,3,\varphi)$ with endvertices $x,y$ under the coloring $\varphi$. We next show that $zu\in E(P)$. Otherwise, we apply a Kempe change on the $(2,3)$-chain containing $zu$. The operation recolors $zu$ with color 3, while keeping the $\varphi(x)$, $|B^\varphi_5|$ and $|S^\varphi_5|$ unchanged. Therefore, we obtain that $xz$ is colored with 1, $zu$ is colored with 3, $d_G(u)=4$, and the set of colors assigned to $u$ is $[4]$, contradicting Subclaim~\ref{subclaim34}. Thus $zu\in E(P)$ as desired. By Subclaim~\ref{subclaim35}, $3\in\varphi(p)$ and $2\in\overline{\varphi}(p)$. Therefore, $p$ is an endvertex of a $(2,3)$-chain. Denote the path $P_p(2,3,\varphi)$ by $P_1$. Since $p\ne x,y$ and the endvertices of $P$ are $x,y$, the paths $P$ and $P_1$ are vertex-disjoint. We apply a Kempe change on $P_u(2,5,\varphi)$, yielding a coloring $\varphi'$. In the set of edges colored with $2$ and $3$, this operation deletes $zu$ and adds $zp$. Since $zu\in E(P)$, deleting $zu$ breaks the path $P$. Thus, $P$ is split into a $(2,3)$-chain containing $x$ and a $(2,3)$-chain containing $y$, with $z$ lying on one of them. The new edge $zp$ joins the $(2,3)$-chain containing $z$ to $P_1$. Since $P_1$ is disjoint from $P$, adding $zp$ does not reconnect $x$ and $y$. Hence, we obtain $P_x(2,3,\varphi')\ne P_y(2,3,\varphi')$.
		
Apply a Kempe change on $P_x(2,3,\varphi')$ to obtain a coloring $\varphi''$. Now we have $2\in\overline{\varphi''}(x)\cap\overline{\varphi''}(y)$, so recolor $xy$ with color 2. Let $\varphi'''$ denote the resulting coloring. The only newly created edge of color 5 is $zu$, which is not a non-local 5-edge, since $d_G(z)\ge 5$. Thus $|B^{\varphi'''}_5|=|B^\varphi_5|-1$, contradicting the choice of $\varphi$.
		
\medskip
\noindent\textbf{Case 2.} $N_G(x)\cap N_G(y)\ne\emptyset$.
		
Let $z\in N_G(x)\cap N_G(y)$. We may assume that
$\varphi(xz)=1$ and $\varphi(yz)=3$.
By Subclaim~\ref{subclaim31} and Subclaim~\ref{subclaim32}, we may write $\varphi(zu)=2$, $\varphi(zw)=4$, and $\varphi(zp)=5$.

By Subclaim~\ref{subclaim31}, $d_H(u)=4$, so $d_G(u)\geq 4$ and $[4]\subseteq\varphi(u)$. If $d_G(u)=4$, then $\varphi(u)=[4]$, contradicting Subclaim~\ref{subclaim34} applied to $yz$ and $zu$. Therefore $d_G(u)\ge 5$, so $u\in V_{\ge5}(G)$. Similarly, $w\in V_{\ge5}(G)$. Since $z\in V_{\ge5}(G)$, the restriction $\Delta(G[V_{\ge5}(G)])\le2$ yields $d_G(p)\le4$.
Since $\varphi(zp)=5$, at least one color in $[4]$ is missing at $p$.
		
\medskip
\noindent\textbf{Subcase 2.1.} $1\in\overline{\varphi}(p)$.
		
Since $1\in\overline{\varphi}(y)$, we have the path $P_y(1,5,\varphi)=(y,yx,x,xz,z,zp,p)$. Applying a Kempe change on $P_y(1,5,\varphi)$ yields that $xz$ is colored with 5, which is not non-local. Let $\varphi'$ denote the resulting coloring. The non-local edge $xy$ is no longer colored with $5$. Thus $|B^{\varphi'}_5|=|B^\varphi_5|-1$, contradicting the choice of $\varphi$.
		
\medskip
\noindent\textbf{Subcase 2.2.} $3\in\overline{\varphi}(p)$.
		
Since $3\in\overline{\varphi}(x)$, we have the path $P_x(3,5,\varphi)=(x,xy,y,yz,z,zp,p)$. Applying a Kempe change on $P_x(3,5,\varphi)$ yields that $yz$ is colored with 5, which is not non-local. Let $\varphi''$ denote the resulting coloring. The non-local edge $xy$ is no longer colored with $5$. Thus $|B^{\varphi''}_5|=|B^\varphi_5|-1$, contradicting the choice of $\varphi$.

Having excluded Subcases 2.1 and 2.2, we may now assume that $1,3\in\varphi(p)$. Since $p$ misses at least one color in $[4]$, it misses color 2 or color 4.

\medskip
\noindent\textbf{Subcase 2.3.} $2\in\overline{\varphi}(p)$.
		
By Subclaim~\ref{subclaim33}, the path $P_x(2,3,\varphi)$ passes through $yz$, hence contains $z$. Since $p$ sees color 3 and misses color 2, $p$ is an endvertex of another $(2,3)$-chain. The two chains are vertex-disjoint.
Applying a Kempe change on $P_p(2,3,\varphi)$ to obtain a new coloring $\varphi'''$ makes $3\in \overline{\varphi'''}(p)$, preserves the relevant edges at $x,y,z$, and preserves all edges of color 5.
Since $3\in \overline{\varphi'''}(x)$, we have the path $P_x(3,5,\varphi''')=(x,xy,y,yz,z,zp,p)$. Applying a Kempe change on $P_x(3,5,\varphi''')$ gives $xy:5\to3$, $yz:3\to5$, and $zp:5\to3$. The only new edge of color 5 is $yz$, which is not non-local because $z\in V_{\ge5}(G)$. The edge $xy$ loses color 5, so that $|B^\varphi_5|$ decreases, contradicting the choice of $\varphi$.
		
\medskip
\noindent\textbf{Subcase 2.4.} $4\in\overline{\varphi}(p)$.
		
By Subclaim~\ref{subclaim33}, the path $P_y(1,4,\varphi)$ passes through $xz$, hence contains $z$. Since $p$ sees color 1 and misses color 4, $p$ is an endvertex of another $(1,4)$-chain. The two chains are vertex-disjoint.
Applying a Kempe change on $P_p(1,4,\varphi)$ to obtain a new coloring $\varphi''''$ makes $1\in \overline{\varphi''''}(p)$, preserves the relevant edges at $x, y, z$, and preserves all edges of color 5.
Since $1\in \overline{\varphi''''}(y)$, we have the path $P_y(1,5,\varphi'''')=(y,yx,x,xz,z,zp,p)$. Applying a Kempe change on $P_y(1,5,\varphi'''')$ gives $xy:5\to1$, $xz:1\to5$, and $zp:5\to1$. The only new edge of color 5 is $xz$, which is not non-local because $z\in V_{\ge5}(G)$. The edge $xy$ loses color 5, so that $|B^\varphi_5|$ decreases, contradicting the choice of $\varphi$.

Since contradictions arise in both Case 1 and Case 2, it follows that $|B^\varphi_5|=0$.
\end{proof}

By the same argument as in the first paragraph of the proof of Theorem \ref{T2}(2), where only the edges colored at most 4 are recolored, we may also assume that $|B^\varphi_t|=0$ for $1\le t\le 4$, and consequently $|B^\varphi_\Delta|=|B^\varphi_{\Delta-1}|=\cdots=|B^\varphi_1|=0$.
Therefore, we can choose a proper $\Delta$-edge-coloring of $G$ that is also locally bounded, and so $\Delta(G)=\chi'(G)\le\chi'_{loc}(G)\le\Delta(G)$, which implies $\chi'_{loc}(G)=\chi'(G)$ as desired.
This completes the proof of Theorem \ref{T3}.
\end{proof}
		
\noindent {\bf Acknowledgements}

The authors are very grateful to Professors Guantao Chen and Gregory Gutin for their valuable suggestions.
This work was supported by Hebei Natural Science Foundation (A2026208005) and Foundation of China Scholarship Council (202508130088).

\end{document}